\documentclass[abstracton,a4paper]{scrartcl}

\usepackage[headsepline]{scrpage2}
\usepackage[latin1]{inputenc}
\usepackage{amsfonts}
\usepackage{amsmath}
\usepackage{amssymb}
\usepackage{amsthm}

\usepackage{url}

\usepackage{graphicx}

\usepackage[pdftex,colorlinks,breaklinks,linkcolor=black,citecolor=black,filecolor=black,menucolor=black,urlcolor=blue,pdfauthor={Lizhen Ji, Peer Kunstmann, and Andreas Weber},pdftitle={}, plainpages=false, pdfpagelabels, bookmarksnumbered=true]{hyperref}

\ihead{L. Ji, P. Kunstmann, A. Weber: Riesz Transformation on Riemannian Manifolds }

\newtheorem{theorem}{Theorem}[section]
\newtheorem{lemma}[theorem]{Lemma}

\newtheorem{corollary}[theorem]{Corollary}

\theoremstyle{definition}

\theoremstyle{remark}

\newcommand{\Hy}{\mathbb{H}}

\newcommand{\Q}{\mathbb{Q}}

\newcommand*\spa{\mathrm{span}}

\newcommand\bP{{\bf P}}

\title{Riesz Transform on Locally Symmetric Spaces and Riemannian Manifolds with a Spectral Gap}
\author{Lizhen Ji\footnote{Email: lji@umich.edu, 
					Address: 1834 East Hall, Ann Arbor, MI 48109-1043, USA.}\\ 
					{\large Department of Mathematics, University of Michigan}
\and Peer Kunstmann\footnote{    Email: peer.kunstmann@math.uni-karlsruhe.de,
						   Address:   Kaiserstr. 89 - 93, 76128 Karlsruhe, Germany.}\\
						   {\large  Institut f\"ur Analysis,
						   Universit\"at Karlsruhe (TH)} 					
\and Andreas Weber\footnote{    Email: andreas.weber@math.uni-karlsruhe.de,
						   Address:   Kaiserstr. 89 - 93, 76128 Karlsruhe, Germany.}\\
						   {\large  Institut f\"ur Algebra und Geometrie,
						   Universit\"at Karlsruhe (TH)} }

\date{}

\begin{document}

\maketitle 
\begin{abstract} 
In this paper we study the Riesz transform on complete and connected Riemannian manifolds 
$M$ with a certain spectral gap in the $L^2$ spectrum of the Laplacian.  
We show that on such manifolds the Riesz transform is $L^p$ bounded for
all $p \in (1,\infty)$.
This generalizes a result by Mandouvalos and Marias and extends a result by Auscher, Coulhon, 
Duong, and Hofmann to the case where zero is an isolated point of the 
$L^2$ spectrum of the Laplacian.\\
%

\end{abstract}		

\section{Introduction and Main Result}

Investigations concerning the question whether the Riesz transform 
$ \nabla \Delta^{-1/2}$ on a general complete Riemannian manifold $M$ is a bounded operator from 
$L^p(M)$ to the space of $L^p$ vector fields, $1< p < \infty,$  have been started with the paper \cite{MR705991} by Strichartz (note that we denote by $\nabla$ the gradient and by $\Delta$ the
Laplace-Beltrami operator on $M$). 
He proved in this paper that the $L^2$ Riesz transform is an isometry for any complete Riemannian manifold $M$ and furthermore, he showed boundedness on $L^p, 1<p<\infty$, if
$M$ is a rank one symmetric space of non-compact type. In general however, the Riesz transform
needs not be bounded for all $p > 1$, see e.g. \cite{MR1458299} where an example for a Riemannian manifold $M$ is given such that the Riesz transform is unbounded on $L^p$ for all
$p>2$.
Note that the $L^p$ boundedness of the Riesz transform on a Riemannian manifold is equivalent
to the existence of a constant $C_p>0$ such that 
$$
 \| |\nabla f|  \|_{L^p} \leq C_p \| \Delta^{1/2}f \|_{L^p}
$$
for all $f\in C_0^{\infty}(M)$ holds.

In \cite{Mandouvalos:2009yq} Mandouvalos and Marias studied the Riesz transform
on certain locally symmetric spaces.  
Inspired by their paper and \cite{MR2119242} we show here how
the $L^p$ boundedness of the Riesz transform 
for a larger class of connected Riemannian manifolds $M$ can be obtained. 
Our only assumptions are a certain spectral gap in the $L^2$ spectrum of the Laplacian and
the following properties of $M$:
\begin{itemize}
 \item[(1.)]  The Riemannian manifold $M$ satisfies the {\em exponential growth property}: 
 For all $r_0>0, x\in M, \theta >1,$ and  $0< r < r_0$  the volume $V(x,r)$ of a ball with center $x$ and
radius $r$ satisfies
\begin{equation}\label{exponential growth property}
V(x,\theta r) \leq Ce^{c\theta}V(x,r)
\end{equation}
for some constants $C\geq 0, c>0$ that only depend on $r_0$. Note that this condition implies
the local volume doubling property, cf.  \cite{MR2119242}.
\item[(2.)] If $p(t,x,y)$ denotes the heat kernel of the heat semigroup $e^{-t\Delta}: L^2(M) \to L^2(M)$
we assume that there exists $C > 0$ such that for all $x,y \in M$ and all $t \in (0,1)$ we have
\begin{equation} \label{DUB}
 p(t,x,x) \leq \frac{C}{V(x,\sqrt{t})}
\end{equation}
and
\begin{equation}\label{gradient}
 \left|\nabla_x p(t,x,y)\right| \leq  \frac{C}{\sqrt{t}\,V(y,\sqrt{t})}.
\end{equation}
\end{itemize}
More precisely, we prove the following theorem.
\begin{theorem}\label{main theorem}
 Let $M$ denote a complete and connected Riemannian manifold satisfying properties
 (\ref{exponential growth property}), (\ref{DUB}), and (\ref{gradient}) from above, 
 and assume that there is a constant $a>0$ such that the following condition for the $L^2$
 spectrum $\sigma(\Delta)$ of $\Delta$ holds:
 $$
   \sigma(\Delta) \subset \{0\} \cup [a,\infty).
 $$
 Then there exist for any $p\in (1,\infty)$  constants $c_p,C_p >0$ 
 such that 
 \begin{equation}\label{inequality Riesz}
   c_p \| \Delta^{1/2}f \|_{L^p} \leq  \| |\nabla f|  \|_{L^p} \leq C_p \| \Delta^{1/2}f \|_{L^p}
 \end{equation}
 for any $f\in C_0^{\infty}(M)$.
\end{theorem}
This result extends \cite[Theorem 1.9]{MR2119242}. The new feature is that we allow
 also $0\in \sigma(\Delta)$. Note that  similar ideas as in our proof  have been used 
 in \cite[Theorem 1.3]{MR1458299}.
 Important examples for manifolds that satisfy the conditions in Theorem \ref{main theorem}  will be given in Section \ref{section proofs}. Note, that properties (\ref{exponential growth property}), (\ref{DUB}), and (\ref{gradient}) are always satisfied if the Ricci curvature of $M$ is bounded from below (cf. also the discussion in  \cite[Section 1.3]{MR2119242}) and hence, these conditions are fulfilled if $M$ is a locally symmetric space. 
\begin{corollary}
 Let $M$ denote a complete and connected Riemannian manifold whose Ricci curvature 
 is bounded from below, and assume that there is a constant $a>0$ with
 $$
   \sigma(\Delta) \subset \{0\} \cup [a,\infty).
 $$
 Then there exist for any $p\in (1,\infty)$  constants $c_p,C_p >0$ 
 such that 
 \begin{equation}
   c_p \| \Delta^{1/2}f \|_{L^p} \leq  \| |\nabla f|  \|_{L^p} \leq C_p \| \Delta^{1/2}f \|_{L^p}
 \end{equation}
 holds for any $f\in C_0^{\infty}(M)$.
\end{corollary}
While Mandouvalos and Marias  in \cite{Mandouvalos:2009yq} proved boundedness of the $L^p$ 
Riesz transform only for $p\in (p_1,p_2)$ with certain $1<p_1<2<p_2<\infty$, we allow here
$p\in (1,\infty)$. 
A standing assumption in the paper \cite{Mandouvalos:2009yq} is that the $L^2$ spectrum
of the Laplacian on a non-compact locally symmetric space $M=\Gamma\backslash X$ 
is of the form $\{\lambda_0,\ldots, \lambda_r\} \cup [||\rho||^2,\infty)$, where 
$\lambda_i < ||\rho||^2, i=0,\ldots,r,$ 
are eigenvalues of finite multiplicity and $\rho$ denotes half the sum of
the positive roots. Note that this assumption needs not be satisfied for all non-compact locally
symmetric spaces. If e.g. the universal covering $X$ is a higher rank symmetric space,
the absolutely continuous part of the $L^2$ spectrum of a non-compact finite volume quotient 
$M=\Gamma\backslash X$ is in many cases of the  form $[a,\infty)$ where $0<a< ||\rho||^2$,
\cite{Ji:2007fk}.
However, it is known that the above cited condition from \cite{Mandouvalos:2009yq} holds true for all non-compact, geometrically finite hyperbolic manifolds $\Gamma\backslash \Hy^n$, \cite{MR661875}.\\
In a first step towards a proof of the boundedness of the Riesz transform the authors show 
in  \cite{Mandouvalos:2009yq} that
any $L^2$ eigenfunction $\varphi_{\lambda_i}$ with respect to some eigenvalue 
$\lambda_i < ||\rho||^2$ is contained in $L^p(M)$ for $p$ in an interval $(p_1,p_2)$ where
\begin{eqnarray*}
  p_1 & = & 2 \left( 1 + \left(1 - \frac{\lambda_i}{|| \rho ||^2} \right)^{1/2} \right)^{-1},\\ 
  p_2 & = & 2  \left( 1 + \left(1 - \frac{\lambda_i}{|| \rho ||^2} \right)^{1/2} \right),\\ 
\end{eqnarray*}
cf. \cite[Theorem 1]{Mandouvalos:2009yq}. Note that the authors use $\lambda_r$ instead of 
$\lambda_i \leq \lambda_r$ in the statement of their result. The proof however shows, that 
$p_1$ and $p_2$ from above can be used. Note further, that $p_1' \geq p_2$, where $p_1'$ 
denotes the conjugate of $p_1$, i.e. $\frac{1}{p_1} + \frac{1}{p_1'} = 1$.
In the proof of this result the authors state that the volume of balls
grows exponentially (formula 2.11) with respect to the radius. 
Unfortunately, this assumption is often not satisfied if $M$ is a locally symmetric space with non-trivial fundamental group. In the case of finite volume it is certainly false. And even
if the volume is infinite the volume growth may be linear (the following example
was communicated to us by Richard D Canary): 
if $M$ is a compact hyperbolic $3$-manifold that fibers over the circle we may 'unwrap' this 
manifold along the circle such that only the fundamental group coming from the fibers survives.
The resulting hyperbolic manifold $\tilde{M}$ is a covering of $M$ that has an infinite cyclic subgroup
of isometries with compact quotient and hence linear volume growth. The existence of hyperbolic
$3$-manifolds that fiber over the circle follows from Thurston's work in \cite{Thurston:yq}
and, of course, from Perelman's  proof of Thurston's geometrization conjecture using 
Ricci flow. However, a detailed look at the proof of  \cite[Theorem 1]{Mandouvalos:2009yq}
shows that only an exponential upper bound on the volume growth is needed and this condition 
is always satisfied for complete locally symmetric spaces $M$.\\
In our proof of the boundedness of the Riesz transform we do not need to make any of the above
assumptions and we will not need specific information about the $L^2$ eigenfunctions.
Nevertheless, we want to mention briefly that in some cases it is possible to enlarge the 
interval $(p_1,p_2)$:
If $M$ is a locally symmetric space with $\Q$-rank one it had been shown in \cite[Theorem 4.1]{Ji:2007fk} that  $p_1=1$ and 
$p_2 = 2( 1 + (1 - \frac{\lambda_i}{|| \rho_{\bP} ||^2})^{1/2})^{-1}$ can be chosen.\\ 
If the injectivity radius of $M$ is strictly positive, it follows immediately from
Taylor's work on the $L^p$ spectrum of the Laplacian \cite[Proposition 3.3]{MR1016445}
that the interval can be enlarged to the right hand side, more precisely, we can choose 
$p_2 = p_1' = 2 ( 1 - (1 - \frac{\lambda_i}{|| \rho ||^2} )^{1/2} )^{-1}$.
Note that the condition on the injectivity radius implies in particular that $M$ has infinite volume.\\ 

\section{Proof and Examples}\label{section proofs}
The following theorem due to Auscher, Coulhon, Duong, and Hofmann is  a main ingredient in our proof of Theorem \ref{main theorem}. 
\begin{theorem}\textup{(\cite[Theorem 1.7]{MR2119242}).} \label{ACDH}
Let $M$ denote a complete Riemannian manifold satisfying the properties (\ref{exponential growth property}), (\ref{DUB}), and (\ref{gradient}), and let $p\in (1,\infty)$. 
Then there exists a constant $C_p > 0$ such that the inequality
\[
 \| | \nabla f | \|_{L^p} + \| f \|_{L^p} \leq 
 C_p \Big( \| \Delta^{1/2}f \|_{L^p} + \| f \|_{L^p} \Big).
\]
holds for all $f\in C_0^{\infty}(M)$.
\end{theorem}
The following lemma can basically be found in the proof of \cite[Theorem 1.3]{MR1458299}.
For the sake of completeness and because we need a variant of the argument below, 
we recall its short proof.
\begin{lemma} \label{lemma fractional}
Assume that $0\notin \sigma(\Delta)$.
Then for any $\alpha \in (0,1)$ the operator $\Delta^{-\alpha}$ defines a  bounded operator on 
$L^p(M), p\in (1,\infty)$. 
\end{lemma}
\begin{proof}
Since $\Delta$ is a selfadjoint operator on $L^2(M)$ and  $\sigma(\Delta)\subset [a,\infty)$
for some $a>0$ it follows from the spectral theorem
$\| e^{-t\Delta} \|_{L^2\to L^2} \leq e^{-at}$.
Furthermore, the semigroup $e^{-t\Delta}: L^2(M)\to L^2(M)$ is positive and 
$L^{\infty}(M)$ contractive. Hence,  this semigroup extends to a 
strongly continuous contraction semigroup
$e^{-t\Delta}: L^p(M)\to L^p(M)$ for any $p\in [1,\infty)$. By interpolation and duality it follows
that for any $p\in [1,\infty)$ we have
$\| e^{-t\Delta} \|_{L^p\to L^p} \leq e^{- 2 \min\{1/p,1/p'\} a t }$ and hence, the integral
\[
 \Delta^{-\alpha}  = \frac{1}{\Gamma(\alpha)}\int_0^{\infty} t^{1-\alpha} e^{-t\Delta} dt.
\]
defines for all $\alpha\in (0,1)$ a bounded operator on $L^p(M), p\in (1,\infty)$. 
\end{proof}
Before we proceed with the proof of Theorem \ref{main theorem}, 
we recall the following general
result that has been proved in \cite{MR705991}. Note that we denote
by $-\Delta_p$ the generator of the heat semigroup $T_p(t) = e^{-t\Delta}$ on $L^p(M)$.
\begin{lemma}\label{core}
For any $1<p<\infty$ and any complete Riemannian manifold $M$, the space
$C^\infty_0(M)$ is a core for the operator $\Delta_p$ in $L^p(M)$, i.e. 
$C^\infty_0(M)$ is dense in ${\cal D}(\Delta_p)$ with respect to the graph norm.
In particular, $C^\infty_0(M)$ is a core for the operator $\Delta_p^{1/2}$.
\end{lemma}
\begin{proof}
Letting $L_p$ denote the $L^p(M)$ closure of $\Delta$, defined on
$C^\infty_0(M)$,  it is shown in \cite{MR705991},
beginning of the proof of Theorem 3.5, that $-L_p$ is dissipative and
generates a contraction semigroup in $L^p(M)$. This semigroup
coincides with the heat semigroup, which implies $L_p=\Delta_p$.
\end{proof}
\paragraph*{Proof of Theorem \ref{main theorem}.}
It is well known (cf. \cite{MR889472,MR990472}) that by duality the first inequality 
in (\ref{inequality Riesz}) is implied by the second and hence we will concentrate 
on the second.\\
The case
$0\notin \sigma(\Delta)$ was already
proved in  \cite[Theorem 1.9]{MR2119242} (cf. also \cite[Theorem 1.3]{MR1458299}): 
By Lemma \ref{lemma fractional} the inverse square root
$\Delta^{-1/2}$  defines a bounded operator on $L^p(M)$ for any $p\in (1,\infty)$ and hence, the 
claim follows from Theorem \ref{ACDH} together with
$\| f \|_{L^p} \leq C \| \Delta^{1/2}f \|_{L^p}$
 for $f\in C_0^{\infty}(M)$.\\
Assume now  $0\in\sigma(\Delta)$. Then $0$ is an isolated point in the spectrum of a 
 self-adjoint operator and hence an eigenvalue. If $f\in L^2(M)$ is a corresponding eigenfunction 
 we may conclude  that $f=const.$ as
 $0=\langle \Delta f,f \rangle = \langle\nabla f, \nabla f \rangle$
 and thus $vol(M) < \infty$.
From H\"older's inequality it follows in particular $L^p(M) \hookrightarrow L^q(M)$
if $1\leq q\leq p\leq \infty$.\\
Let us define for $p\in [1,\infty)$
\[
 L^p_0(M) = \left\{ f\in L^p(M) :  \langle 1, f\rangle=\int_M f dx = 0 \right\}. 
\]
Then $L^p_0(M)$ is an invariant subspace for the semigroup 
$T_p(t)=e^{-t\Delta}: L^p(M)\to L^p(M)$ and $L^p(M) = L^p_0(M) \oplus \spa\{ 1\}$.
Furthermore, if we denote by $-\Delta_0$ the generator of $e^{-t\Delta}: L^2_0(M) \to L^2_0(M)$
we have $\Delta_0 = \Delta\big|_{{\cal D}(\Delta_2)\cap L^2_0(M)}$.
As $\sigma(\Delta_0) \subset [a,\infty)$ for some $a>0$ we can prove as in Lemma \ref{lemma fractional}
that $\Delta_0^{-\alpha}$ is bounded on $L^p_0(M)$ for all $\alpha\in (0,1)$ and $p\in (1,\infty)$.
Note that we need here to assure that the spaces $L^p_0(M), p\in [1,\infty),$ interpolate
correctly, i.e. $[L^{p_1}_0(M), L^{p_2}_0(M)]_{\theta} =  L^{p_\theta}_0(M)$ (see \cite{MR1328645}), where
$\frac{1}{p_{\theta}} = \frac{1-\theta}{p_1} + \frac{\theta}{p_2}$. This however follows from the fact
that 
\[
  P: L^p(M) \to L^p(M), \,  f\mapsto f - \frac{1}{vol(M)} \int_M f dx
\]
is a continuous projection onto $L^p_0(M)$ together with \cite[1.2.4]{MR1328645}.
If $M$ is compact then inequality (\ref{inequality Riesz}) follows from Theorem \ref{ACDH} 
since $Pf\in C^\infty_0(M)\cap L^p_0(M)$ and (\ref{inequality Riesz}) holds for the
constant function $g:=f-Pf$.\\
If $M$ is non-compact but of finite volume then $Pf\in C^\infty_0(M)$
is no longer true in general. 
But by Lemma \ref{core}, Theorem \ref{ACDH} holds by
approximation for all $g\in {\cal D}(\Delta_p^{1/2})$ with the same
constant. 
In particular, it holds for all $g\in C^\infty_0(M)+\spa\{1\}$. Now we can finish the proof as
in the compact case since $f\in C^\infty_0(M)$ implies 
$Pf\in \left(C^\infty_0(M)+\spa\{1\}\right)\cap L^p_0(M)$.

%
%
\subsection*{Examples}
The conditions in Theorem \ref{main theorem} are satisfied for the following
manifolds.
\paragraph*{Compact manifolds.}
If $M$ is compact, its $L^2$ spectrum is discrete and $0\in \sigma(\Delta)$.
\paragraph*{Locally symmetric spaces.}
Let  $M= \Gamma\backslash G/K$  denote a locally symmetric space whose universal covering
is a symmetric space of non-compact type.\\
If the critical exponent 
 $\delta(\Gamma)$ satisfies $\delta(\Gamma) < 2||\rho||$ ($\rho$ denotes half the sum of the positive roots) then the bottom of the $L^2$
 spectrum is bounded from below by a strictly positive constant \cite{MR2019974,Weber:2007fk}. 
 Note that in this case $vol(M)=\infty$.\\
 If $M$ is non-compact with finite volume, the $L^2$ spectrum equals 
 $\sigma(\Delta)= \{0, \lambda_1,\ldots, \lambda_r\} \cup [b,\infty)$ for some $b>0$, see e.g. 
 M\"uller's article in \cite{:2008xq}.\\
 Important examples in this context are all non-compact hyperbolic manifolds
 $M=\Gamma\backslash \Hy^n$ where $\Gamma$ is geometrically finite. Then, 
 in the case of finite and infinite volume of $M$, the 
 $L^2$ spectrum is of the form
 $\{\lambda_1,\ldots, \lambda_r\} \cup \left[\frac{(n-1)^2}{4},\infty\right)$, cf. \cite{MR661875}.
 \paragraph*{Manifolds with cusps of rank one.} These non-compact manifolds with 
 finite volume satisfy also the condition  
  $\sigma(\Delta)= \{0, \lambda_1,\ldots, \lambda_r\} \cup [b,\infty)$ for some $b>0$,
  cf. \cite{MR891654}.
 \paragraph*{Homogeneous Spaces.}
 Consider $M = \Gamma\backslash G$ where G is a semisimple Lie group endowed 
 with an invariant metric and $ \Gamma \subset G$ is a non-uniform arithmetic lattice in $G$. 
 Then, the spectrum equals 
 $\sigma(\Delta)= \{0, \lambda_1,\ldots, \lambda_r\} \cup [b,\infty)$ for some $b>0$, cf.
 \cite{MR701563}.\\

\noindent Note that the continuous $L^2$ spectrum of the Laplacian on a complete Riemannian manifold is invariant under compact perturbations of the Riemannian metric.
This follows from the so-called decomposition principle in \cite{MR544241}. Hence, the spectral
gap condition in the above mentioned cases is still satisfied after compact perturbations.
\paragraph*{Is zero in the spectrum?}
 
 This question was raised in various contexts. Let us recall the following results
 which immediately give further examples.
 \begin{theorem}[\cite{MR683635}]
   Let $M$ have infinite volume. Suppose that there is a constant $\kappa\geq 0$
   such that $Ricci_M \geq -\kappa^2$. Then $0\notin \sigma(\Delta)$ if and only if
   $M$ is open at infinity.
 \end{theorem} 
 Recall that a manifold with infinite volume is called open at infinity if and only if the Cheeger constant
 is strictly positive.
 \begin{theorem}[\cite{MR656213}]
   Let $X$ be a normal covering of a compact manifold $M$ with covering group $\Gamma$.
   Then $0\in \sigma(\Delta)$ if and only if $\Gamma$ is amenable.    
 \end{theorem}
 \begin{theorem}[\cite{MR1250269} or \cite{MR945013}]
  Let $M$ denote a Riemannian manifold with Ricci curvature bounded from below and 
  with empty cut locus. If there is a point $x_0\in M$ such that the volume form 
  $\sqrt{g(x_0,\zeta)}$
  of the  Riemannian metric grows exponentially in every direction (with respect to geodesic normal
  coordinates $(x,\zeta)$), the bottom of the $L^2$ spectrum is strictly positive.
 \end{theorem}
 
\paragraph*{Acknowledgements} We want to thank Michel Marias for valuable remarks on a first version of this paper. We also want to thank the referee for many suggestions which led to a great improvement of this paper. Lizhen Ji was partially supported by NSF grant DMS 0604878.

\bibliographystyle{amsplain}
\bibliography{dissertation,hypercyclic,symmetricSpaces}

\end{document}